\documentclass[graybox]{svmult}

\usepackage{mathptmx}       
\usepackage{helvet}         
\usepackage{courier}        
\usepackage{type1cm}        
\usepackage{graphicx,makeidx,color}
\usepackage{epsfig}
\usepackage{array}
\usepackage{booktabs}
\usepackage{multicol}        
\usepackage[bottom]{footmisc}

\usepackage{amssymb}

\newtheorem{de}{Definition}[section]

\newcommand{\W}{\bf W}

\newcommand{\uu}{\mathrm{u}}
\newcommand{\x}{\mathrm{x}}
\newcommand{\y}{\mathrm{y}}
\newcommand{\z}{\mathrm{z}}
\newcommand{\s}{\mathrm{s}}
\newcommand{\EE}{\mathbf{E}}

\newcommand{\DD}{\mathbf{D}}
\newcommand{\dd}{\mathrm{d}}
\newcommand{\R}{\mathbf{R}}
\newcommand{\be}{\begin{equation}}
\newcommand{\ee}{\end{equation}}

\makeindex             

\begin{document}

\title*{Multidimensional Sensitivity Analysis of Large-scale Mathematical Models}
\titlerunning{Multidimensional Sensitivity Analysis}
\author{Ivan Dimov and Rayna Georgieva}
\institute{Ivan Dimov and Rayna Georgieva \at Department of Parallel Algorithms, IICT, Bulgarian Academy of Sciences, Acad. G. Bonchev 25 A, 1113 Sofia, Bulgaria, \email{ivdimov@bas.bg, rayna@parallel.bas.bg}}

%
%

\maketitle
\abstract{Sensitivity analysis (SA) is a procedure for studying  how sensitive are the output results of large-scale mathematical models to some uncertainties of the input data.  The models are described as a system of partial differential equations. Often such systems contain a large number of input parameters. Obviously, it is important to know how sensitive is the solution to some uncontrolled variations or uncertainties in the input parameters of the model. Algorithms based on analysis of variances technique (ANOVA) for calculating numerical indicators of sensitivity and computationally efficient Monte Carlo integration techniques have recently been developed by the authors. They have been successfully applied to sensitivity studies of air pollution levels calculated by the Unified Danish Eulerian Model (UNI-DEM) with respect to several important input parameters.
In this paper a comprehensive theoretical and experimental study of the Monte Carlo algorithm based on \textit{symmetrised shaking} of Sobol sequences has been done. It has been proven that this algorithm  has an optimal rate of convergence for functions with continuous and bounded second derivatives in terms of probability and mean square error. Extensive numerical experiments with Monte Carlo, quasi-Monte Carlo (QMC) and scrambled quasi-Monte Carlo algorithms based on Sobol sequences are performed to support the theoretical studies and to analyze applicability of the algorithms to various classes of problems. The numerical tests show that the Monte Carlo algorithm based on \textit{symmetrised shaking} of Sobol sequences gives reliable results for multidimensional integration problems under consideration.}

\section{Introduction}
\label{Intr}

Most existing methods for providing SA rely
on special assumptions connected to the behavior of the model
(such as linearity, monotonicity and additivity of the
relationship between model input and model output) \cite{STCR04}.
Such assumptions are often applicable to a large range of mathematical
models. At the same time there are models that include significant nonlinearities and/or
stiffness. For such models assumptions about linearity and additivity are not
applicable. This is especially true when one deals with non-linear systems of partial differential equations. The numerical study and results reported in this paper have been done by using a large-scale mathematical model called Unified Danish Eulerian Model (UNI-DEM) \cite{ZD06,ZDG96}. The model enables us to study the transport of air pollutants and other species over a large geographical region. The system of partial differential equations describes the main physical processes, such as advection, diffusion, deposition, as well as chemical and photochemical processes between the studied species. The emissions, and the quickly changing meteorological conditions are also described. The non-linearity of the equations are mainly introduced when modeling chemical reactions \cite{ZD06}. If the model results are sensitive to a given process, one can describe it mathematically in a more adequate way, or more precisely. Thus, the goal of our study is to increase the reliability of the results produced by the model, and to identify processes that must be studied more carefully, as well as to find input parameters that need to be measured with a higher precision.  A careful sensitivity analysis is needed in order to decide where and how simplifications of the model can be
made. That's why it is important to develop and study more adequate and reliable methods for sensitivity analysis. A good candidate for reliable sensitivity analysis of models containing nonlinearity is the variance based method \cite{STCR04}.
The idea of this approach is to estimate
how the variation of an input parameter or a group of inputs contributes into
the variance of the model output.  
As a measure of this analysis we use the \textit{total sensitivity indices (TSI)}  
(see, Section \ref{Setting}) described as multidimensional integrals:

\be \label{int}
I = \int_{\Omega} g(\x)p(\x) \, \dd \x, \,\,\, \Omega \subset \R^d,
\ee
where $g(\x)$ is a square integrable function in $\Omega$ and
$p(\x) \geq 0$ is a probability density function, such that
$\int_{\Omega} p(\x) \, \dd \x = 1$.

That's why it is important to deal with efficient numerical methods for high-dimensional integration. The progress in the area of sensitivity analysis is closely connected to the progress in reliable algorithms for multidimensional integration.

\section{Problem Setting}
\label{Setting}

\subsection{Modeling and Sensitivity}
\label{mod-sens}

Assume that the mathematical model can be presented as a
 function
\begin{equation}\label{math_model}
\uu = f(\x), \quad \mbox{where} \quad \x = (x_1,x_2,\ldots,x_d) \in
U^d \equiv [0;1]^d
\end{equation}
is the vector of input parameters with a joint {\bf p}robability
{\bf d}ensity {\bf f}unction (p.d.f.) $p(\x) = p(x_1,\ldots,x_d)$.
Assume also that the input variables are independent (non-correlated)
and the density function $p(\x)$ is known, even if $x_i$ are not actually random
variables (r.v.).
The \textit{total sensitivity index} \cite{HS96} provides a
measure of the total effect of a given parameter, including all
the possible coupling terms between that parameter and all the
others.
The {\bf t}otal {\bf s}ensitivity {\bf i}ndex (TSI) of an
input parameter $x_i, i \! \in \! \{1,\ldots,d\}$ is defined in the
following way \cite{HS96,Sobol01}:
\be\label{TSI} S_{i}^{tot} =
S_i + \displaystyle \sum_{l_1 \ne i} S_{i l_1} + \sum_{l_1, l_2
\ne i, l_1 < l_2} S_{i l_1 l_2} + \ldots + S_{i l_1 \ldots
l_{d-1}},
\ee
where $S_i$ is called \textit{the main effect
(first-order sensitivity index)} of $x_i$ and $S_{i l_1 \ldots
l_{j-1}}$ is the  $j^{\mbox{-th}}$ order sensitivity index.
The higher-order terms describe the
interaction effects between the unknown input parameters
$x_{i_1},\ldots,x_{i_{\nu}}, \nu \in \{2,\ldots,d\}$ on the
output variance.

The method of global SA used in this work  is based on a decomposition
of an integrable model function $f$ in the $d$-dimensional factor
space into terms of increasing dimensionality \cite{Sobol01}: \be\label{mod-f}
f(\x) = f_0 + \sum_{{\nu}=1}^d \sum_{l_1 < \ldots < l_{\nu}}
f_{l_1 \ldots l_{\nu}}(x_{l_1}, x_{l_2}, \ldots, x_{l_{\nu}}), \ee
where $f_0$ is a constant.
The representation (\ref{mod-f}) is
referred to as the ANOVA-representation of the model function $f(\x)$ if each term is chosen to
satisfy the following condition \cite{Sobol01}:
$$
\int_0^1 f_{l_1 \ldots l_{\nu}}(x_{l_1},x_{l_2},\ldots,x_{l_{\nu}}) \dd
x_{l_k} = 0, \quad 1 \leq k \leq {\nu}, \quad {\nu}=1,\dots,d.
$$
Let us mention the fact that if the whole presentation (\ref{mod-f}) of
the right-hand site is used, then it doesn't simplify the problem. The hope is
that a truncated sequence $f_0 + \sum_{{\nu}=1}^{d_{tr}} \sum_{l_1 < \ldots < l_{\nu}}
f_{l_1 \ldots l_{\nu}}(x_{l_1}, x_{l_2}, \ldots, x_{l_{\nu}})$,
where $d_{tr}<d$ (or even $d_{tr}<<d$), can be considered as a good approximation to the
model function $f$.

The quantities
\be\label{variances}
\DD = \displaystyle \int_{U^d} f^2(\x) \dd \x - f_0^2, \quad
\DD_{l_1 \ \ldots \ l_{\nu}} = \int f^2_{l_1 \ \ldots \
l_{\nu}} \dd x_{l_1} \ldots \dd x_{l_{\nu}}
\ee
are the so-called total and partial variances respectively and are
obtained after squaring and integrating over $U^d$ the equality
(\ref{mod-f}) on the assumption that $f(\x)$ is a square integrable function
(thus all terms in (\ref{mod-f}) are also square integrable functions).
Therefore, the total variance of the model output is split
into partial variances in the analogous way as the
model function, that is the unique ANOVA-decomposition:
$
\DD = \sum_{{\nu}=1}^d \sum_{l_1 < \ldots < l_{\nu}} \DD_{l_1 \ldots l_{\nu}}.
$
The use of probability theory concepts is based on the assumption that the input
parameters are random variables distributed in
$U^d$ that defines $f_{l_1 \ \ldots \ l_{\nu}}(x_{l_1}, x_{l_2}, \ldots, x_{l_{\nu}})$
also as random variables with variances (\ref{variances}).
For example $f_{l_1}$ is presented by a conditional expectation:
$
f_{l_1}(x_{l_1}) = \EE (\uu|x_{l_1}) - f_0$ and respectively
$\DD_{l_1} = \DD[f_{l_1}(x_{l_1})] = \DD[\EE (\uu|x_{l_1})].
$
Based on these  assumptions about the model function and the output variance,
the following quantities
\be \label{sen_in}
S_{l_1 \ \ldots \ l_{\nu}} =
\displaystyle \frac{\DD_{l_1 \ \ldots \ l_{\nu}}}{\DD}, \quad {\nu} \in
\{1,\ldots, d\}
\ee
are referred to as the global sensitivity indices
\cite{Sobol01}.
Based on the formulas (\ref{variances})-(\ref{sen_in}) it is clear that the
mathematical treatment of the problem of providing global
sensitivity analysis consists in evaluating total sensitivity
indices (\ref{TSI}) of corresponding order that, in turn, leads to
computing multidimensional integrals of the form (\ref{int}).
It means that to obtain $S_{i}^{tot}$ in general,
 one needs to compute $2^d$ integrals of type (\ref{variances}).
As we discussed earlier the basic assumption underlying representation
(\ref{mod-f}) is that the basic features of the model functions
(\ref{math_model}) describing typical real-life problems can be
presented by low-order subsets of input variables,
containing terms of the order up to $d_{tr}$, where $d_{tr}<d$ (or even $d_{tr}<<d$). Therefore,
based on this assumption, one can assume that the dimension of the initial problem can
be reduced.

The procedure for computing global sensitivity indices (see \cite{Sobol01}) is based
on the following representation of the variance
\be \label{var_y}
\DD_{\y} :\displaystyle \DD_{\y} = \int \ f(\x) \ f(\y,\z') \dd \x \dd \z' - f_0^2,
\ee
where $\y = (x_{k_1},\ldots,x_{k_m}), \ 1 \leq k_1 < \ldots < k_m \leq
d, \ $ is an arbitrary set of $m$ variables ($1 \leq m \leq d-1$) and $\z$ is the set of $d - m$ complementary variables, i.e.  $\x
= (\y,\z)$. The  equality (\ref{var_y}) enables the construction of a
Monte Carlo algorithm for evaluating $f_0, \DD$ and $\DD_{\y}$:
$$
\begin{array}{ll}
\displaystyle \frac{1}{n}\ \sum_{j=1}^n \ f(\xi_j)
\stackrel{P}{\longrightarrow} f_0, & \qquad
\displaystyle \frac{1}{n}\ \sum_{j=1}^n \ f(\xi_j) \ f(\eta_j,\zeta'_j)\stackrel{P}{\longrightarrow} \DD_{\y} + f_0^2,  \\
\displaystyle \frac{1}{n}\ \sum_{j=1}^n \ f^2(\xi_j)
\stackrel{P}{\longrightarrow} \DD + f_0^2, & \qquad \displaystyle
\frac{1}{n}\ \sum_{j=1}^n \ f(\xi_j) \
f(\eta'_j,\zeta_j)\stackrel{P}{\longrightarrow} \DD_{\z} + f_0^2,
\end{array}
$$
where $\xi = (\eta,\zeta)$ is a random sample and $\eta$ corresponds to the input subset denoted by $\y$.

Instead of randomized (Monte Carlo) algorithms for computing the above
sensitivity parameters one can use deterministic quasi-Monte Carlo algorithms,
or randomized quasi-Monte Carlo \cite{LLT08,LL02}. Randomized (Monte Carlo)
algorithms have proven to be very efficient in solving
multidimensional integrals in composite domains \cite{2008D,So73}. At the same time
the QMC based on well-distributed Sobol sequences can be
considered as a good alternative to Monte Carlo algorithms, especially
for smooth integrands and not very high \textit{effective dimensions} (up to $d=15$) \cite{Ku11}.
Sobol $\Lambda \Pi_{\tau}$ are good candidates for efficient
QMC algorithms. Algorithms based on $\Lambda \Pi_{\tau}$ sequences while being deterministic,
mimic the pseudo-random sequences used in Monte Carlo integration.
One of the problems with $\Lambda \Pi_{\tau}$ sequences is that they
may have \textit{bad} two-dimensional projection. In this context \textit{bad} means that
the distribution of the points is far away from the uniformity.
If such projections are used in a certain computational problem, then the lack of uniformity may provoke a substantial lost
of accuracy. To overcome this problem randomized
QMC can be used. There are several ways of randomization and \textit{scrambling} is
one of them.
The original motivation of scrambling \cite{HH03,Owen95} aims toward obtaining more
uniformity for quasi-random sequences in high dimensions, which can be checked
via two-dimensional projections. Another way of randomisation is to \textit{shake} the quasi-random
points according to some procedure.
Actually, the scrambled algorithms  obtained by \textit{shaking} the quasi-random
points can be considered as Monte Carlo algorithms with a special choice of the density function. It's
a matter of definition. Thus, there is a reason to be able to compare two classes of algorithms:
\textit{deterministic} and \textit{randomized}.

\section{Complexity in Classes of Algorithms}

One may pose the task to consider and compare
two classes of algorithms:
\textit{deterministic algorithms} 
and
\textit{randomized (Monte Carlo) algorithms}. 
Let $I$ be the desired value of the integral. Assume for a given
r.v. $\theta$ one can prove that the mathematical
expectation satisfies $\EE \theta = I$. Suppose that the mean
value of $n$ values of $\theta$: $\theta^{(i)}, \ i=1,\dots,n$
is considered as a Monte Carlo approximation to the solution:
$
\bar{\theta}_n = 1/n\sum_{i=1}^n \theta^{(i)} \approx I,
$
where $\theta^{(i)} (i = 1, 2, \dots, n)$ correspond to values (realizations) of a r.v. $\theta$.
In general, a certain randomized algorithm can produce the
result with a given probability error.
So, dealing with randomized algorithms one has to accept that the
result of the computation can be true only with a certain (although
high) probability. In most practical computations it is
reasonable to accept an error estimate with a probability smaller
than $1$.

Consider the following integration problem:

\be \label{int-pr} S(f):=I = \int_{U^d} f(\x) d\x, \ee where
$\x \equiv (x_1, \dots , x_d) \in U^d \subset \R^d$ and
$f \in C(U^d)$ is an integrable function on $U^d$. The computational
problem can be considered as a mapping of function $f: \{[0,1]^d
\rightarrow \R \}$ to $ \R$: $ S(f):  f \rightarrow
\R, $ where $ S (f) = \int_{U^d} f(\x) d\x$ and $f\in F_0 \subset
C(U^d)$. We refer to $S$ as the solution operator. The elements of
$F_0$ are the data, for which the problem has to be solved; and for
$f \in F_0, \hskip 0.3cm S (f) $ is the exact solution. For a given
$f$, we want to compute exactly or approximately $S (f)$. One may be interested
in cases when the integrand $f$ has a higher regularity. It is because
in many cases of practical computations $f$ is smooth and has high order
bounded derivatives. If this is the case, then is it reasonable to try to exploit
such a smoothness. To be able to do that we need to define the functional class
$F_0 \equiv {\W}^k(\| f\|; {U}^d)$ in the following way:

\begin{de}
Let $d$ and $k$ be
integers,  $d, k\geq 1$.  We consider the class ${\W}^k(\| f\|;
{U}^d)$ (sometimes abbreviated to ${\W}^k$) of real functions $f$
defined over the unit cube $U^d=[0,1)^d$, possessing all the
partial derivatives
$
\frac{\partial ^rf(\x)}{\partial x_1^{\alpha _1}\dots \partial x_d^{\alpha _d}
}, \,\,\, \alpha _1 + \dots + \alpha _d=r\leq k ,
$
which are continuous when $r<k$ and bounded in $\sup $ norm when
$r=k.$ The semi-norm $\left\| \cdot \right\|$ on ${\W}^k$ is defined
as
$$
\left\| f\right\| =\sup \left\{ \left| \frac{\partial ^kf(\x)}{\partial
x_1^{\alpha _1}\dots \partial x_d^{\alpha _d}}\right|,\,\,\,\, \alpha _1+\dots
+\alpha _d=k,\,\,\, \x\equiv (x_1,...,x_d)\in U ^d\right\} .
$$
\end{de}
We call a {\it quadrature formula} any expression of the form
$$ A^D(f,n) =
\sum_{i=1}^n c_i f(\x^{(i)}),
$$
which approximates the value of the
integral $S (f)$. The real numbers $c_i \in \R$ are called weights
and the $d$ dimensional points $\x^{(i)} \in U^d$ are called nodes. It is
clear that for fixed weights $c_i$ and nodes $\x^{(i)}\equiv (x_{i,1},\ldots , x_{i,d})$ the quadrature
formula $A^D(f,n)$ may be used to define an algorithm with an integration error $
err(f,A^D) \equiv \int_{U^d}f(\x)d\x- A^D(f,n)
$. We call a
{\it randomized quadrature formula} any formula of the following
kind:
$
A^R(f,n) = \sum_{i=1}^n \sigma_i f(\xi^{(i)}),
$
where $\sigma_i$
and $\xi^{(i)}$ are random weights and nodes respectively. The algorithm $A^R(f,n)$
belongs to the class of randomized (Monte Carlo) denoted by $\cal A^R$.

\begin{de} \label{pr-err-def}
Given a randomized (Monte Carlo) integration formula for the functions from the space ${\W}^k$ we define the integration error
$$
err(f,A^R) \equiv \int_{U^d}f(\x)d\x- A^R(f,n)
$$
by the probability error $\varepsilon_P (f)$ in the sense that $\varepsilon_P (f)$ is the least possible real number, such that
$$
Pr\left( \left| err(f,A^R)\right| <\varepsilon_P (f)\right) \geq P,
$$
and the mean square error
$$
r(f)=\left\{ E \left[ err^2(f,A^R)\right]\right\} ^{1/2}.
$$
\end{de}

We assume that it suffices to obtain an $\varepsilon_P (f)$-approximation
to the solution with a probability $0<P<1$. If we allow equality, i.e., $0<P \le 1$ in
Definition \ref{pr-err-def}, then  $\varepsilon_P (f)$ can be used as an accuracy measure
for both randomized and deterministic algorithms. In such a way it is consistent to
consider a wider class $\cal A$ of algorithms that contains both classes:
randomized and deterministic algorithms.

\begin{de}\ \ Consider the set $\cal A$ of algorithms $A$:
$$ {\cal A} = \{A:Pr(|err(f,A)| \le \varepsilon) \ge c\}, \ \ \ A\in \{A^D, A^R\}, \ \ 0 < c < 1 $$ that solve a given
problem with an integration error $err(f,A)$.
\end{de}

In such a setting it is correct to compare randomized algorithms with algorithms
based on low discrepancy sequences like Sobol $\Lambda \Pi_{\tau}$ sequences.

\section{The Algorithms}

The algorithms we study are based on Sobol $\Lambda \Pi_{\tau}$ sequences.

\subsection{$\Lambda \Pi_{\tau}$ Sobol Sequences}
\label{Sob-seq}

$\Lambda \Pi_{\tau}$ sequences are \textit{uniformly
distributed sequences} (u.d.s.) The term \textit{u.d.s.} was introduced
by Hermann Weyl in 1916 \cite{We16}.
For practical purposes an u.d.s.  must be found that satisfied three
requirements \cite{So73,So89}:
(i) the best asymptote as $n \rightarrow \infty$,
(ii) well distributed points for small $n$, and
(iii)  a computationally inexpensive algorithm.

All  $\Lambda \Pi_{\tau}$-sequences given in
\cite{So89} satisfy the first requirement. Suitable distributions such as $\Lambda \Pi_{\tau}$ sequences are also called
$(t,m,s)$-nets and $(t,s)$-sequences in base $b \ge 2$. To introduce them,
define first an elementary $s$-interval in base $b$ as a subset of ${U}^s$ of the form
$ E= \prod_{j=1}^s \left [ \frac{a_j}{b^{d_j}}, \frac{a_j+1}{b^{d_j}}\right],
$ where $a_j, d_j\ge 0$ are integers and $a_j< b^{d_j}$ for all $j\in \{1, ...,s\}$.
Given two integers $0\le t \le m$, a $(t,m,s)$-net in base $b$ is a sequence $\x^{(i)}$ of $b^m$ points of ${U}^s$
such that $Card \ E \cap \{\x^{(1)}, \ldots , \x^{(b^m)}\}=b^t$ for any elementary interval $E$ in base $b$ of hypervolume $\lambda(E) = b^{t-m}$.
Given a non-negative integer $t$, a $(t,s)$-sequence in base $b$ is an infinite sequence
of points $\x^{(i)}$ such that for all integers $k \ge 0, m \ge t$, the sequence
$\{\x^{(kb^m)}, \ldots, \x^{((k+1){b^m-1})}\}$ is a $(t,m,s)$-net in base $b$.

I. M. Sobol \cite{So73} defines his $\Pi_{\tau}$-meshes and $\Lambda \Pi_{\tau}$ sequences, which are $(t,m,s)$-nets
and $(t,s)$-sequences in base $2$ respectively. The terms $(t,m,s)$-nets and $(t,s)$-sequences
in base $b$ (also called Niederreiter sequences) were introduced in 1988 by H. Niederreiter \cite{Ni88}.

To generate the $j$-th component of the points in a Sobol sequence, we need to choose a
primitive polynomial of some degree $s_j$ over the Galois field of two elements GF(2)
$ P_j = x^{s_j} + a_{1,j} x^{s_j -1} + a_{2,j} x^{s_j -2}+ \ldots + a_{s_j-1,j} x +1,$
where the coefficients $a_{1,j},\ldots, a_{s_j-1,j}$ are either $0$ or $1$.
A sequence of positive integers $\{m_{1,j}, m_{2,j},\ldots\}$
are defined by the recurrence relation
$$m_{k,j}=2a_{1,j}m_{k-1,j}\oplus 2^2 a_{2,j} m_{k-2,j}\oplus \dots \oplus 2^{s_j} m_{k-s_j,j}\oplus m_{k-s_j,j},$$
where $\oplus$ is the bit-by-bit \textit{exclusive-or} operator.
The values $m_{1,j}, \dots , m_{s_j,j}$ can be chosen freely provided that each $m_{k,j}, 1 \le k \le s_j$,
is odd and less than $2^k$. Therefore, it is possible to construct different Sobol sequences
for the fixed dimension $s$. In practice, these numbers must be chosen very carefully to obtain really efficient Sobol sequence generators \cite{SAKK11}.
The so-called direction numbers $\{v_{1,j}, v_{2,j},\dots \}$ are defined by $\displaystyle v_{k,j}= \frac{m_{k,j}}{2^k}$.
Then the $j$-th component of the $i$-th point in a Sobol sequence, is given by
$x_{i,j} = i_1 v_{1,j}\oplus i_2 v_{2,j} \oplus \ldots,$
where $i_k$ is the $k$-th binary digit of $i = (\dots i_3i_2i_1)_2$.
Subroutines to compute these points can be found in \cite{BF88,So79}. The work \cite{LMRS88} contains more details.

\subsection{The Monte Carlo Algorithms based on Modified Sobol Sequences - MCA-MSS}

One of the algorithms based on a procedure of \textit{shaking} was proposed recently in \cite{2010DG}.
The idea is that we take a Sobol $\Lambda \Pi_{\tau}$
point (vector) $\x$ of dimension $d$. Then $\x$ is considered
as a centrum of a sphere with a radius $\rho$. A random point $\xi \in U^d$ uniformly distributed on the sphere
is taken. Consider a random variable $\theta$ defined as a value of the integrand at that random point, i.e., $\theta = f(\xi)$.
Consider random points $\xi^{(i)}(\rho) \in U^d, i=1, \ldots, n$. Assume $\xi^{(i)}(\rho) = \x^{(i)} + \rho \omega^{(i)}$,
where $\omega^{(i)}$ is a unique uniformly distributed vector in $U^d$. The radius $\rho$ is
relatively small $\rho <<  \frac{1}{2^{d_j}}$, such that $\xi^{(i)}(\rho)$
is still in the same elementary $i^{th}$ interval
$E_i^d = \prod_{j=1}^d \left [ \frac{a_j^{(i)}}{2^{d_j}}, \frac{a_j^{(i)}+1}{2^{d_j}}\right]
$, where the pattern $\Lambda \Pi_{\tau}$ point $\x^{(i)}$ is. We use a subscript $i$ in $E_i^d$ to
indicate that the $i$-th $\Lambda \Pi_{\tau}$ point $\x^{(i)}$ is in it. So, we assume that if
$\x^{(i)} \in E_i^d$, then $\xi^{(i)}(\rho) \in E_i^d$ too.

It was proven in \cite{2010DG} that
the mathematical expectation of the random variable $\theta = f(\xi)$ is equal to the value of the integral
(\ref{int-pr}), that is
$ \EE \theta = S(f) = \int_{U^d} f(\x) d\x.$
This result allows for defining a randomized algorithm. One can take the Sobol
$\Lambda \Pi_{\tau}$ point $\x^{(i)}$ and \textit{shake} it somewhat. \textit{Shaking} means to
define random points $\xi^{(i)}(\rho) = \x^{(i)} + \rho \omega^{(i)}$ according to the procedure described above.
For simplicity the algorithm described above is called MCA-MSS-1.

The probability error of the algorithm MCA-MSS-1 was analysed in \cite{12DGOZ-CAM}.
It was proved that for integrands with continuous and bounded first derivatives, i.e. $f \in F_0 \equiv {\W}^1(L; {U}^d)$,
where $L=\| f\|$, it holds
$$
err(f,d)\leq c_{d}^{^{\prime }}\left\| f\right\|
n^{^{-\frac 12-\frac 1d}} \quad \mbox{and} \quad
r(f,d)\leq c_{d}^{^{\prime \prime }}\left\| f\right\|
n^{^{-\frac 12-\frac 1d}},
$$
where the constants $c_{d}^{^{\prime }}$ and $c_{d}^{^{\prime \prime }}$
do not depend on $n$.

In this work a modification of algorithm MCA-MSS-1 is proposed and analysed.
The new algorithm will be called MCA-MSS-2.

It is assumed that $n=m^d$,  $m \geq 1$. The unit cube
$U^d$ is divided into $m^d$ disjoint sub-domains, such that they coincide with the elementary
$d$-dimensional subintervals
defined in Subsection \ref{Sob-seq}
$
U^d=\bigcup_{j=1}^{m^d}K_j,\,\,\, {\rm{ where }}\,\,\,K_j=
\prod_{i=1}^d[a_i^{(j)},b_i^{(j)}),
$
with $b_i^{(j)}-a_i^{(j)}=\displaystyle \frac {1}{m}$ for all $i=1, \dots ,d$.

In such a way in each $d$-dimensional sub-domain $K_j$ there is exactly one
$\Lambda \Pi_{\tau}$ point $\x^{(j)}$. Assuming that after shaking, the random point
stays inside $K_j$, i.e., $\xi^{(j)}(\rho) = \x^{(j)} + \rho \omega^{(j)}\in K_j$ one
may try to exploit the smoothness of the integrand in case if the integrand has second continuators
and bounded derivatives, i.e., $f \in F_0 \equiv {\W}^2(L; {U}^d)$.

Then, if $p(x)$ is a probability density function, such that $\int_{U^d} p(\x) d\x =1$, then
$$
\int_{K_j} p(\x) d\x = p_j \le \frac{c_1^{(j)}}{n},
$$
where $c_1^{(j)}$ are constants.
 If $d_j$ is the diameter of $K_j$, then
$$
d_j= \sup_{x_1, x_2 \in K_j}|x_1 - x_2|\le \frac{c_2^{(j)}}{n^{1/d}},
$$
where $c_2^{(j)}$ are another constants.

In the particular case when the subintervals are with edge $1/m$ for all constants we have:
$c_1^{(j)}=1$ and $c_2^{(j)}=\sqrt{d}$. In each sub-domain $K_j$ the central point is denoted by $\s^{(j)}$, where $\s^{(j)}= (s_1^{(j)}, s_2^{(j)}, \dots , s_d^{(j)})$.

Suppose two random points $\xi^{(j)}$ and $\xi^{(j)'}$ are chosen, such that $\xi^{(j)}$ is selected during our procedure used in MCA-MSS-1.
The second point $\xi^{(j)'}$ is chosen to be symmetric to $\xi^{(j)}$ according to the central point $s^{(j)}$ in each cube
$K_j$. In such away the number of random points is $2 m^d$. One may
calculate all function values $f(\xi^{(j)})$ and $f(\xi^{(j)'})$, for $j=1, \dots, m^d$ and
approximate the value of the integral in the following way:
\be \label{QF2}
I(f)\approx \frac {1}{2m^d} \sum_{j=1}^{2n} \left[ f(\xi^{(j)}) + f(\xi^{(j)'}) \right].
\ee
This estimate corresponds to MCA-MSS-2.
Later on it will be proven that this algorithm
has an optimal rate of convergence for functions with second bounded derivatives, i.e., for functions
$f \in F_0 \equiv {\W}^2(L; {U}^d)$, while the algorithm MCA-MSS-1 has an optimal rate of convergence
for functions with first bounded derivatives: $f \in F_0 \equiv {\W}^1(L; {U}^d)$.

One can prove the following

\begin{theorem}
The quadrature formula (\ref{QF2}) constructed above for integrands $f$ from ${\W}^2(L; {U}^d)$ satisfies
$$
err(f,d)\leq \widetilde{c}_{d}^{\ \prime }\left\| f\right\|
n^{^{-\frac 12-\frac 2d}}
$$
and
$$
r(f,d)\leq \widetilde{c}_{d}^{\ \prime \prime }\left\| f\right\|
n^{^{-\frac 12-\frac 2d}},
$$
where the constants $\widetilde{c}_{d}^{\ \prime }$ and $\widetilde{c}_{d}^{\ \prime \prime }$
do not depend on $n$.
\end{theorem}

\begin{proof}
One can see that
$$
\EE\left\{ \frac {1}{2m^d}\sum^{2n}_{j=1} \left[ f(\xi^{(j)})+ f(\xi^{(j)'})\right] \right\} =\int_{{U}^d}f(\x)d\x.
$$

For the fixed $\Lambda \Pi_{\tau}$ point $\x^{(j)}\in K_j$ one can use the $d$-dimensional Taylor formula to present the function $f(\x^{(j)})$ in $K_j$ around the central point  $\s^{(j)}$.
For simplicity the superscript of the argument $(j)$ will be omitted assuming that the formulas are written for the $j^{th}$ cube $K_j$:
$$ 
f(\x) =  \sum_{n_1=0}^{\infty} \sum_{n_2=0}^{\infty} \dots \sum_{n_d=0}^{\infty}
\frac{(x_1 - s_1)^{n_1} \dots (x_d - s_d)^{n_d}}{n_1! \dots n_d!}
 \frac{\partial^{n_1+\dots+n_d} f}{\partial {x_1}^{n_1}\dots \partial {x_d}^{n_d}} (s_1, \dots, s_d).
$$
Now, one can write this formula  at previously defined random points $\xi$ and $\xi'$ both belonging to $K_j$. In such a way we have:

\be \label{xi}
f(\xi) = f(\s) + [Df(\s)]^T \ (\xi - \s) + \frac{1}{2!}(\xi -\s)^T \ D^2f(\s) \ (\xi - \s) + \ldots,
\ee

\be \label{xi'}
f(\xi') = f(\s) + [Df(\s)]^T \ (\xi' - \s) + \frac{1}{2!}(\xi' -\s)^T \ D^2f(\s) \ (\xi' - \s) + \ldots,
\ee
where $Df(\s)$ is the gradient of $f$ evaluated at $\x=\s$ and $D^2f(\s)$ is the Hessian matrix, i.e.,
$$
D^2f(\s) =  \left[ \begin{array}{cccc} \frac{\partial^2 f}{\partial x_1^2} & \frac{\partial^2 f}{\partial x_1\partial x_2}
 & \ldots & \frac{\partial^2 f}{\partial x_1\partial x_d} \\
 \frac{\partial^2 f}{\partial x_2\partial x_1} & \frac{\partial^2 f}{\partial x_2^2}
 & \ldots & \frac{\partial^2 f}{\partial x_2\partial x_d} \\
 \vdots & \vdots & \ddots & \vdots  \\
\frac{\partial^2 f}{\partial x_d\partial x_1} & \frac{\partial^2 f}{\partial x_d\partial x_2}
 & \ldots & \frac{\partial^2 f}{\partial x_d^2}\end{array} \right].
$$
Summarising (\ref{xi}) and (\ref{xi'}) one can get
$$
f(\xi) + f(\xi') = 2f(\s) + \frac{1}{2!} \ D^2f(\s) \ \left( (\xi -\s)^T (\xi - \s) + (\xi' -\s)^T (\xi' - \s)\right) + \dots.
$$
Because of the symmetry there is no member depending on the gradient $Df(\s)$ in the previous formula. If we consider
the variance $ \DD [f(\xi) + f(\xi')]$ taking into account that the variance of the constant $2f(\s)$ is zero, then we will get
$$
\begin{array}{ll}
\DD [f(\xi) + f(\xi')] & =   \DD \left[ \frac{1}{2} D^2f(\s) \left( (\xi -\s)^T (\xi - \s) + (\xi' -\s)^T (\xi' - \s)\right) + \dots \right] \\
& \le \EE \left[   \frac{1}{2} D^2f(\s) \left( (\xi -\s)^T (\xi - \s) + (\xi' -\s)^T (\xi' - \s)\right) + \dots  \right] ^2.
\end{array}
$$
Now we return back to the notation with superscript taking into account that the above consideration is just for an arbitrary sub-domain $K_j$. Since $f\in {\W}^2(L; U^d)$, $\|f\|\le L_j$ and $L=\|f\|$ is the majorant for all
$L_j$, i.e., $L_j \le L$ for $j=1, \ldots, n$. Obviously,
there is a point ($d$-dimensional vector) $\eta \in K_j$, such that
$$
f(\xi) + f(\xi') \le 2f(\s) + \frac{1}{2} \ D^2f(\s) \ \left( (\eta -\s)^T (\eta - \s) +
(\eta -\s)^T (\eta - \s)\right)
$$
and the variance can be estimated from above in the following way:
\begin{eqnarray*}
\DD [f(\xi) + f(\xi')] & \le & L_j^2 \sup_{x_1^{(j)},x_2^{(j)}}\left|x_1^{(j)}- x_2^{(j)}\right|^4  \le  L_j^2(c_2^{(j)})^4 n^{-4/d}.
\end{eqnarray*}
Now the variance of $\theta_n=\sum_{j=1}^n \theta^{(j)}$ can be estimated:
\begin{eqnarray}\label{varian}
\DD\theta_n = \sum_{j=1}^n p_j^2 \DD [f(\xi) + f(\xi')] & \le & \sum_{j=1}^n (c_1^{(j)})^2 n^{-2} L_j^2(c_2^{(j)})^4 n^{-4/d} \nonumber  \\
                                                  & \le & \left( L_j c_1^{(j)} c_2^{(j)2}\right)^2 n^{-1-4/d}.
\end{eqnarray}
Therefore $ \qquad
r(f,d)\leq \widetilde{c}_{d}^{\ \prime \prime } \left\| f\right\|
n^{^{-\frac 12-\frac 2d}}.
$
The application of the Tchebychev's
inequality to the variance (\ref{varian}) yields
$$
\varepsilon (f,d)\leq \widetilde{c}_{d}^{\ \prime }\left\|
f\right\| n^{^{-\frac 12-\frac 2d}}
$$
for the probable error $\varepsilon$,
where $\widetilde{c}_{d}^{\ \prime }=\sqrt{2d}$, which concludes the proof.

\end{proof}

One can see that the Monte Carlo algorithm  MCA-MSS-2 has an optimal rate of
convergence for functions with continuous and bounded second derivative \cite{2008D}. This means that the rate of convergence
($\displaystyle n^{-\frac 12 - \frac{2}{d}}$) can not be improved for the functional class ${\W}^2$ in the class of the randomized algorithms $\cal A^R$.

Note that both MCA-MSS-1 and MCA-MSS-2 have one control parameter, that is the radius $\rho$ of the \textit{sphere of shaking}. At the same time, to be able to efficiently use this control parameter
one should increase the computational complexity. The problem is that after \textit{shaking} the random point may leave the multidimensional sub-domain. That's why
after each such a procedure one should be checking if the random point is still in the same sub-domain. It is clear that the procedure of checking if a random point is inside the
given domain is a computationally  expensive procedure when one has a large number of points. A small modification of MCA-MSS-2 algorithm allows to overcome this difficulty. If we
just generate a random point $\xi^{(j)} \in K_j$  uniformly distributed inside $K_j$ and after that take the symmetric point $\xi^{(j)'}$ according to the central point $s^{(j)}$, then this procedure
will simulate the algorithm  MCA-MSS-2. Such a completely  randomized approach simulates algorithm  MCA-MSS-2, but the \textit{shaking} is with different radiuses $\rho$ in each sub-domain.
We will call this algorithm MCA-MSS-2-S, because this approach looks like the stratified symmetrised Monte Carlo. Obviously, MCA-MSS-2-S is less expensive than MCA-MSS-2, but there is not
such a control parameter like the radius $\rho$, which can be considered as a parameter randomly chosen in each sub-domain  $K_j$.

It is important to notice that all three algorithms MCA-MSS-1, MCA-MSS-2 and MCA-MSS-2-S have optimal (unimprovable) rate of convergence for the corresponding functional classes, that is MCA-MSS-1 is
optimal in $F_0 \equiv {\W}^1(L; {U}^d)$ and both MCA-MSS-2 and MCA-MSS-2-S are optimal in $F_0 \equiv {\W}^2(L; {U}^d)$.

We also will be considering the known Owen Nested Scrambling Algorithm \cite{Owen95} for which it is
proved that the rate of convergence is $n^{-3/2} (log\ n)^{(d-1)/2}$, which is very good but still not optimal even for integrands in $F_0 \equiv {\W}^1(L; {U}^d)$. One can see that if the \textit{logarithmic function} from the estimate can be omitted, then the rate will became optimal. Let us mention that it is still not proven that the above estimate is exact, that is, we do not know if the  \textit{logarithm} can be omitted. It should be mentioned that the proved convergence rate for the Owen Nested Scrambling Algorithm
improves significantly the rate for the unscrambled nets, which is $n^{-1} (log\ n)^{d-1}$.
That's why it is important to compare numerically our algorithms MCA-MSS with the Owen Nested Scrambling.
The idea of Owen nested scrambling is based on randomization of a single digit at each iteration.
Let $x^{(i)} = (x_{i,1}, x_{i,2}, \ldots, x_{i,s}), \ i=1, \ldots, n$ be  quasi-random numbers in $[0, 1)^s$, and let $z^{(i)} =(z_{i,1}, z_{i,2}, \ldots, z_{i,s})$
be the scrambled version of the point $x^{(i)}$. Suppose that each $x_{i,j}$ can be
represented in base $b$ as $x_{i,j} = (0.x_{i1,j} \ x_{i2,j} \ldots x_{iK,j} \ldots)_b$ with $K$ being the number of digits to be scrambled.
Then nested scrambling proposed by Owen \cite{Owen95,Owen02} can be defined as follows:
$z_{i1,j} = \pi_{\bullet}(x_{i1,j})$, and $z_{il,j} = \pi_{\bullet x_{i1,j} x_{i2,j} \ldots x_{i {l-1},j}}
(x_{il,j})$, with independent permutations $\pi_{\bullet x_{i1,j} x_{i2,j} \ldots x_{i {l-1},j}}$
for $l \geq 2$. Of course, $(t,m,s)$-net remains $(t, m, s)$-net under nested scrambling. However,
nested scrambling requires $b^{l-1}$ permutations to scramble the $l$-th digit. Owen scrambling
(nested scrambling), which can be applied to all $(t, s)$-sequences, is powerful; however,
from the implementation point-of-view, nested scrambling or so-called path dependent
permutations requires a considerable amount of bookkeeping, and leads to more problematic
implementation. There are various versions of scrambling methods based on digital permutation, and
the differences among those methods are based on the definitions of the $\pi_l$'s. These
include Owen nested scrambling \cite{Owen95,Owen02}, Tezuka's generalized Faure sequences \cite{Tez95}, and Matousek's linear
scrambling \cite{Mat98}.

\section{Case-study: Variance-based Sensitivity Analysis of the Unified Danish Eulerian Model}

The input data for the sensitivity analysis performed in this paper
has been obtained during runs of
a large-scale mathematical model for remote transport of air pollutants
(\textbf{Uni}fied \textbf{D}anish \textbf{E}ulerian \textbf{M}odel,
UNI-DEM, \cite{ZD06}). The model enables us to study
concentration variations in time of a high number of air pollutants and
other species over a large geographical region (4800 $\times$ 4800 km),
covering the whole of Europe, the Mediterranean and some parts of Asia
and Africa. Such studies are important for environmental protection, agriculture,
health care. The model presented as a system of partial differential equations 
describes the main processes in the atmosphere including
photochemical processes between the studied species, the emissions, the
quickly changing meteorological conditions.
Both non-linearity and stiffness of the equations are mainly introduced
when modeling chemical reactions \cite{ZD06}. The chemical scheme used in the
model is the well-known condensed CBM-IV (Carbon Bond Mechanism). Thus,
the motivation to choose UNI-DEM is that
it is one of the models of atmospheric chemistry, where the chemical
processes are taken into account in a very accurate way.

This large and complex task is not suitable for direct numerical
treatment.
For the purpose of numerical solution
it is split into submodels, which represent the main physical and
chemical processes. The sequential splitting \cite{M85}
is used in the production version of the model, although other splitting
methods have also been considered and
implemented in some experimental versions \cite{DFHZ04,DOZ05}. Spatial
and time discretization makes each of the above submodels
a huge computational task, challenging for the most powerful
supercomputers available nowadays. That is why
parallelization has always been a key point in the computer
implementation of DEM since its very early stages.

Our main aim here is to study the sensitivity of the ozone concentration
according to the rate variation of some chemical reactions. We consider
the chemical rates
to be the input parameters and the concentrations of
pollutants to be the output parameters.

\section{Numerical Results and Discussion}

Some numerical experiments are performed to study experimentally various properties of the algorithms.
The expectations based on theoretical results are that for non-smooth
functions MCA-MSS algorithms based on the \textit{shaking} procedures outperform the QMC even for relatively low dimensions. It is also interesting to observe how behave the randomized QMC based on scrambled Sobol sequences.

For our numerical tests we use the following non-smooth integrand:

\be\label{non-smooth-f}
f_1(x_1,x_2,x_3,x_4) = \displaystyle \sum_{i=1}^{4} | (x_i-0.8)^{-1/3} |,
\ee
for which even the first derivative does not exist. Such kind of applications appear also in some important problems in financial mathematics. The referent value of the integral $S(f_1)$ is approximately equal to $7.22261$.

To make a comparison we also consider an integral with a smooth integrand:
\be \label{smooth-f}
f_2(x_1,x_2,x_3,x_4) =  x_1 \ x_2^2 \ {\bf e}^{x_1  x_2} \sin x_3 \cos x_4.
\ee
The second integrand (\ref{smooth-f}) is an infinitely smooth function with a referent value of the integral $S(f_2)$ approximately equal to $0.10897$. The integration domain in both cases is $U^4 = [0,1]^4$.

Some results from the numerical integration tests with a smooth (\ref{smooth-f}) and a non-smooth  (\ref{non-smooth-f}) integrand are presented in Tables \ref{tab_smooth-1} and \ref{tab_non-smooth} respectively. As a measure of the efficiency of the algorithms both the relative error (defined as the absolute error divided by the referent value) and computational time are shown.
\begin{table}[t]\caption{Relative error and computational time for numerical integration of a smooth function ($S(f_2) \approx 0.10897$).}\label{tab_smooth-1}
\begin{tabular}{p{0.63cm}|p{1.097cm}p{1.097cm}|p{1.097cm}p{1.097cm}|p{1.097cm}p{1.097cm}|p{1.097cm}p{1.097cm}p{1.097cm}}\hline
$n~~$        & \multicolumn{2}{c|}{SFMT}      & \multicolumn{2}{c|}{Sobol QMCA}   & \multicolumn{2}{c|}{Owen scrambling} & \multicolumn{3}{c}{MCA-MSS-1}   \\ \cline{2-3} \cline{4-5} \cline{6-7} \cline{8-10}
 &  Rel.  & Time  &  Rel.  & Time  & Rel.  & Time  &  $\rho$        &  Rel.  & Time  \\
 &  error & (s)   &  error & (s)   & error & (s)   &  $\times 10^3$ &  error &  (s)  \\ \hline
$10^2$    & 0.0562      & 0.002    & 0.0365      & $< 0.001$ & 0.0280     & 0.001 & 3.9  & 0.0363 &  0.001 \\ 
          &             &          &             &           &            &       & 13   & 0.0036 &  0.001 \\ 
$10^3$    & 0.0244      & 0.004    & 0.0023      & $ 0.001$  & 0.0016     & 0.001 & 1.9  & 0.0038 &  0.010 \\ 
          &             &          &             &           &            &       & 6.4  & 0.0019 &  0.010 \\ 
$10^4$    & 0.0097      & 0.019    & 0.0009      & 0.002     & 0.0003     & 0.003 & 0.8  & 0.0007 &  0.070 \\ 
          &             &          &             &           &            &       & 2.8  & 0.0006 &  0.065  \\ \hline
\end{tabular}
\end{table}
For generating Sobol quasi-random sequences the algorithm with Gray code implementation \cite{AS79} and sets of direction numbers proposed by Joe and Kuo \cite{JK08} are used. The MCA-MSS-1 algorithm \cite{2010DG} involves generating random points uniformly distributed on a sphere with radius $\rho$. One of the best available random number generators, SIMD-oriented Fast Mersenne Twister (SFMT) \cite{SM08,SFMT} 128-bit pseudo-random number generator of period $2^{19937}-1$  has been used to generate the required random points. SFMT algorithm is a very efficient implementation of the Plain Monte Carlo method \cite{So73}. The radius $\rho$ depends on the integration domain, number of samples and minimal distance between Sobol deterministic points $\delta$. We observed experimentally that the behaviour of the relative error of numerical integration is significantly influenced by the fixed radius of spheres. That is why the values of the radius $\rho$ are presented according to the number of samples $n$ used in our experiments, as well as to a fixed coefficient, \textit{radius coefficient} $\kappa=\rho / \delta$. The latter parameter gives the ratio of the radius to the minimal distance between Sobol points.
\begin{table}[t]\caption{Relative error and computational time for numerical integration of a non-smooth function ($S(f_1) \approx 7.22261$).}\label{tab_non-smooth}
\begin{tabular}{p{0.65cm}|p{1.097cm}p{1.097cm}|p{1.097cm}p{1.097cm}|p{1.097cm}p{1.097cm}|p{1.097cm}p{1.097cm}p{1.097cm}} 
\hline
 $n$   & \multicolumn{2}{c|}{SFMT}      & \multicolumn{2}{c|}{Sobol QMCA} & \multicolumn{2}{c|}{Owen scrambling} & \multicolumn{3}{c}{MCA-MSS-1}   \\  \cline{2-3} \cline{4-5} \cline{6-7} \cline{8-10}
 &  Rel.  & Time  &  Rel.  & Time  & Rel.  & Time  &  $\rho$        &  Rel.  & Time  \\
 &  error & (s)   &  error & (s)   & error & (s)   &  $\times 10^3$ &  error &  (s)  \\  \hline
$10^3$    & 0.0010      & 0.011    & 0.0027      &   0.001   & 0.0021     & 0.002 & 1.9  & 0.0024 &  0.020 \\ 
          &             &          &             &           &            &       & 6.4  & 0.0004 &  0.025 \\ 
$7.10^3$  & 0.0009      & 0.072    & 0.0013      & 0.009     &0.0003      & 0.011 & 1.0  & 0.0004 &  0.110  \\ 
          &             &          &             &           &            &       & 3.4  & 0.0005 &  0.114  \\ 
$3.10^4$  & 0.0005      & 0.304    & 0.0003      & 0.032     & 0.0003     & 0.041 & 0.6  & 0.0001 &  0.440  \\ 
          &             &          &             &           &            &       & 1.9  & 0.0002 &  0.480  \\ 
$5.10^4$  & 0.0007      & 0.513    & 0.0002      & 0.053     &2{\bf e}-05 & 0.066 & 0.4  & 7{\bf e}-05 &  0.775  \\ 
          &             &          &             &           &            &       & 1.4  & 0.0001 &  0.788  \\ \hline
\end{tabular}
\end{table}
The code of scrambled quasi-random sequences used in our studies is taken from the collection of NAG C Library \cite{nag_cl}. This implementation of scrambled quasi-random sequences is based on TOMS Algorithm 823 \cite{HH03}. In the implementation of the scrambling there is a possibility to make a choice of three methods of scrambling: the first is a restricted form of Owen scrambling \cite{Owen95},
the second based on the method of Faure and Tezuka \cite{FT00},
and the last method combines the first two (it is referred to as a \emph{combined} approach).

Random points for the MCA-MSS-1 algorithm have been generated using the original Sobol sequences and modeling a random direction in $d$-dimensional space. The computational time of the calculations with pseudo-random numbers generated by SFMT (see columns labeled as \emph{SFMT} and \emph{MCA-MSS} in Tables \ref{tab_smooth-1} and \ref{tab_non-smooth} has been estimated for all $10$ algorithm runs.

Comparing the results in Tables \ref{tab_smooth-1} and \ref{tab_non-smooth}  one  observes that
\begin{itemize}
\item all algorithms under consideration are efficient and converge with the expected rate of convergence;
\item in the case of \textit{smooth functions}, the Sobol algorithm is better than SFMT (the relative error is up to $10$ times smaller than for SFMT);
\item the scrambled QMC and MCA-MSS-1 are much better than the classical Sobol algorithm; in many cases  even the simplest \textit{shaking} algorithm MCA-MSS-1 gives a higher accuracy than the scrambled algorithm.
\item in case of \textit{non-smooth functions} SFMT algorithm implementing the plain Monte Carlo method is better than the Sobol algorithm for relatively small samples ($n$);
\item in the case of \textit{non-smooth functions} our Monte Carlo \textit{shaking} algorithm MCA-MSS-1 gives similar results as the scrambled QMC; for several values of $n$ we  observe advantages for MCA-MSS-1 in terms of accuracy;
\item both MCA-MSS-1 and scrambled QMC are better than SFMT and Sobol quasi MC algorithm in the case of non-smooth functions.
\end{itemize}

Another observation is that for the chosen integrands  the scrambling algorithm does not outperform the algorithm with the original Sobol points, but the scrambled algorithm and Monte Carlo algorithm MCA-MSS-1 are more stable with respect to relative errors for relatively small values of $n$.

The facts we observed that some further improvements of implementation of more refined \textit{shaking} algorithms MCA-MSS-2 and  MCA-MSS-2-S may be expected for relatively smooth integrands. That's why we
compare Sobol QMCA with MCA-MSS-2 and MCA-MSS-2-S, as well as  with simplest \textit{shaking} algorithm MCA-MSS-1 (see Table \ref{tab_smooth2}).
\begin{table}[t]\caption{Relative error and computational time for
numerical integration of a smooth function ($S(f) \approx 0.10897
$).}\label{tab_smooth2}
{\tabcolsep=3pt \baselineskip=5pt \setlength{\extrarowheight}{0.05cm}
\begin{tabular}{p{1.45cm}|p{0.88cm}p{0.97cm}|p{0.88cm}p{0.88cm}p{0.88cm}|p{0.88cm}p{0.88cm}|p{0.88cm}p{0.88cm}} \hline
\# of points $n$        & \multicolumn{2}{c|}{Sobol QMCA}  &
\multicolumn{3}{c|}{MCA-MSS-1}  & \multicolumn{2}{c|}{MCA-SMS-2}  &
\multicolumn{2}{c}{MCA-SMS-2-S}  \\
\cline{2-3} \cline{4-6}
\cline{7-8} \cline{9-10}
(\# of double                 &  Rel.      & Time        &$\rho$      &  Rel.
& Time      & Rel.        & Time  & Rel.        & Time  \\
      points $2n$)               &  error      & (s)          &$\times 10^3$&  error
&  (s)      & error      & (s)  & error      &  (s)
 \\  \hline
$2^9$              & 0.0059      & $<0.001$    & 2.1        & 0.0064
& 0.009      & 0.0033      & 0.010 & 0.0016      & 0.005 \\
$(2 \times 2^9)$    &            &              & 6.4        & 0.0061
& 0.010      & 0.0032      & 0.010 &            &      \\ 
$2^{10}$            & 0.0035      & 0.002        & 1.9        & 0.0037
& 0.010      & 9{\bf e}-05 & 0.020 & 0.0002      & 0.007 \\
$(2 \times 2^{10})$ &            &              & 6.4        & 0.0048
& 0.010      & 0.0002      & 0.020 &            &      \\ 
$2^{16}$            & 2{\bf e}-05 & 0.027        & 0.4        & 3{\bf
e}-05  & 1.580      & 7{\bf e}-06 & 1.340 & 9{\bf e}-06 &  0.494\\
$(2 \times 2^{16})$ &            &              & 1.2        & 0.0001
& 1.630      & 5{\bf e}-06 & 1.380 &            &      \\ \hline
\end{tabular}
}
\end{table}
The results show that the simplest \textit{shaking} algorithm MCA-MSS-1 gives relative errors similar to
errors of the Sobol QMCA, which is expected since the $\Lambda \Pi_{\tau}$ Sobol sequences are already quite well distributed. That's why one should not expect improvement for a very smooth integrand. But the \textit{symmetrised shaking} algorithm MCA-MSS-2 improves the relative error. The effect of this improvement is based on the fact that the second derivatives of the integrand exists, they are bounded and the construction of the MCA-MSS-2 algorithm gives a better convergence rate of order $O(n^{-1/2-2/d})$. The same convergence rate has the algorithm MCA-MSS-2-S, but the latter one does not allow to control the value of the radius of \textit{shaking}. As expected MCA-MSS-2-S gives better results than MCA-MSS-1.
The relative error obtained by MCA-MSS-2 and MCA-MSS-2-S are of the same magnitude (see Table \ref{tab_smooth2}). The advantage of MCA-MSS-2-S is that its computational complexity is much smaller.
A comparison of the relative error and computational complexity for different values of $n$ is presented
in Table \ref{tab_smooth}. To have a fair comparison we have to consider again a smooth function (\ref{smooth-f}).
\begin{table}[t]\caption{Relative error and computational time for
numerical integration of a smooth function ($S(f) \approx 0.10897
$).}\label{tab_smooth}
\centering
{\tabcolsep=3pt \baselineskip=5pt \setlength{\extrarowheight}{0.05cm}
\begin{tabular}{p{1.5cm}|p{1.2cm}p{1.2cm}|p{1.2cm}p{1.2cm}p{1.2cm}|p{1.2cm}p{1.2cm}}
\hline
$n$        & \multicolumn{2}{c|}{Sobol QMCA} &
\multicolumn{3}{c|}{ MCA-MSS-1 }&
\multicolumn{2}{c}{MCA-MSS-2-S }  \\  \cline{2-3} \cline{4-6} \cline{7-8}
 &  Rel.  & Time  &  $\rho$        &  Rel.  & Time &  Rel.  & Time  \\
 &  $\times 10^3$ &  error &  (s)  &  error & (s)  & error & (s)
 \\  \hline
$2 \times 4^4$    & 0.0076      & $<0.001$    &
2.1  & 0.0079 &  $<0.001$ & 0.0016    & 0.005 \\
(512)            &            &              &   6.4  & 0.0048 &  $<0.001$         &      &
 \\ 
$2 \times 6^4$    & 0.0028      & 0.001      &
1.2  & 0.0046 &  0.030    & 0.0004    & 0.009\\
(2592)            &            &              &       4.1  & 0.0046 &  0.030      &      &
 \\ 
$2 \times 8^4$    & 0.0004      & 0.004        &
0.9  & 0.0008 &  0.090  & 0.0002    & 0.025 \\
(8192)            &            &              &         2.9  & 0.0024 &  0.090    &      &
 \\ 
$2 \times 10^4$  & 0.0002      & 0.008        &
0.6  & 0.0001 &  0.220  &5{\bf e}-05 & 0.070 \\
(20000)          &            &              &   2.0  & 0.0013 &  0.230           &      &
\\
$2 \times 13^4$  & 0.0001      & 0.022        &
0.4  & 0.0001 &  0.630  &4{\bf e}-06 & 0.178 \\
(57122)          &            &              & 1.2  & 0.0007 &  0.640            &      &
 \\
$2 \times 14^4$  & 5{\bf e}-06 & 0.029        &
0.4  & 1{\bf e}-05 &  0.860  &1{\bf e}-05 & 0.237 \\
(76832)          &            &              & 1.2  & 0.0005 &  0.880            &      &
 \\
$2 \times 15^4$  & 8{\bf e}-06 & 0.036         &
0.4  & 0.0001 &  1.220 &9{\bf e}-07 & 0.313 \\
(101250)          &            &              & 1.2  & 0.0005 &  1.250            &      &
  \\ \hline
\end{tabular}
}
\end{table}
The observation is that MCA-MSS-2-S algorithm outperforms the simplest \textit{shaking} algorithm MCA-MSS-1 in terms of relative error and complexity.

After testing the algorithms under consideration on the \textit{smooth} and \textit{non-smooth functions}
we studied the efficiency of the algorithms on \textit{real-life} functions obtained after
running UNI-DEM. Polynomials of 4-th degree with 35 unknown coefficients are used to approximate
 the \textit{mesh functions}
containing the model outputs.

We use various values of the number of points that corresponds to situations when one needs to compute the sensitivity measures with different accuracy.  We have computed results for $g_0$ ($g_0$ is the integral
over the integrand $g(x) = f(x) - c$, $f(x)$ is the approximate model function of UNI-DEM, and  $c$ is a constant obtained as a Monte Carlo estimate of $f_0$, \cite{SobM07}), the total variance $\mathbf{D}$, as well as total sensitivity indices $S_i^{tot}, i=1,2,3$. The above mentioned parameters are presented in Table  \ref{tab6600}. Table \ref{tab6600} presents the results obtained for a relatively low sample size  $n=6600$.

One can notice that for most of the sensitivity parameters the simplest \textit{shaking} algorithm MCA-MSS-1 outperforms the scrambled Sobol sequences, as well as the algorithm based on the $\Lambda \Pi_{\tau}$ Sobol sequences in terms of accuracy. For higher values of sample sizes this effect is even stronger.
\begin{table}[t]\caption{Relative error (in absolute value) and computational time for estimation of sensitivity indices of input parameters using various Monte Carlo and quasi-Monte Carlo approaches ($n = 6600, c \approx 0.51365, \delta \approx 0.08$).}\label{tab6600}
\centering
{\tabcolsep=3pt \baselineskip=100pt
\begin{tabular}{p{1.9cm}|p{2.15cm}|p{2.15cm}|p{2.15cm}p{2.15cm}}
\hline
Estimated   & Sobol QMCA  &  Owen scrambling    & \multicolumn{2}{c}{MCA-MSS-1} \\ \cline{4-5}
quantity    &   &    &  $\rho$ & Rel. error     \\ \hline
$g_0$           &1{\bf e}-05  &  0.0001    & 0.0007 & 0.0001 \\ 
                &         &                & 0.007  & 6{\bf e}-05\\ 
$\mathbf{D}$    & 0.0007  &  0.0013        & 0.0007 & 0.0003        \\ 
                &         &                & 0.007  & 0.0140        \\ 
$S_1^{tot}$     & 0.0036  &  0.0006        & 0.0007 & 0.0009        \\ 
                &         &                & 0.007  & 0.0013        \\ 
$S_2^{tot}$     & 0.0049  &  6{\bf e}-05   & 0.0007 & 2{\bf e}-05   \\ 
                &         &                & 0.007  & 0.0034        \\ 
$S_3^{tot}$     & 0.0259  &  0.0102        & 0.0007 & 0.0099        \\ 
                &         &                & 0.007  & 0.0211        \\ \hline
\end{tabular}
}
\end{table}

One can clearly observe that the simplest \textit{shaking} algorithm MCA-MSS-1 based on modified Sobol sequences improves the error estimates for non-smooth integrands. For smooth functions modified algorithms MCA-MSS-2 and MCA-MSS-2-S give better results than  MCA-MSS-1.
Even for relatively large radiuses $\rho$ the results are good in terms of accuracy. The reason is that centers of spheres are very well uniformly distributed by definition. So that, even for large values of radiuses of \textit{shaking} the generated random points continue to be well distributed.
We should stress on the fact that for relatively low number of points ($< 1000$) the algorithm based on modified Sobol sequences gives results with a high accuracy.

\section{Conclusion}

A comprehensive theoretical and experimental study of the Monte Carlo algorithm MCA-MSS-2 based on \textit{symmetrised shaking} of Sobol sequences has been done. The algorithm combines properties of two of the best available approaches - Sobol quasi-Monte Carlo integration and a high quality SFMT pseudo-random number generator. It has been proven that this algorithm has an optimal rate of convergence for functions with continuous and bounded second derivatives in terms of probability and mean square error.

A comparison with the scrambling approach, as well as with the Sobol quasi-Monte Carlo algorithm and the algorithm using SFMT generator has been provided for numerical integration of smooth and non-smooth integrands. The algorithms mentioned above are tested numerically also for computing sensitivity measures for UNI-DEM model to study sensitivity of ozone concentration according to variation of chemical rates. All algorithms under consideration are efficient and converge with the expected rate of convergence. It is important to notice that the Monte Carlo algorithm MCA-MSS-2 based on modified Sobol sequences when \textit{symmetrised shaking} is used has a unimprovable rate of convergence and gives reliable numerical  results.

\section*{Acknowledgment}
The research reported in this paper is partly supported by the Bulgarian NSF
Grants DTK 02/44/2009 and DMU 03/61/2011.
%
%

\end{document}